\documentclass[a4paper,12pt]{amsart} 
\usepackage{euscript,eufrak,verbatim}
\usepackage[final]{epsfig}
\usepackage{subfigure}
\usepackage{psfrag}

\def\e{\varepsilon}

\def\Z{\mathbb{Z}} \def\R{\mathbb{R}} \def\S{\mathbb{S}}

\newenvironment{proofof}[1]{\noindent {\bf Proof of #1.}}{ \hfill\qed\\ }
\renewenvironment{proof}{\noindent {\bf Proof.}}{ \hfill\qed\\ }

\def\bs{\bar{s}}

\def\t{\hat{t}}
\def\hl{\hat{\lambda}}
\def\hb{\hat{\beta}}
\def\bv{\bar{v}}

\def\v{\vec{v}}
\def\u{\vec{u}}
\def\phi{\varphi}
\def\be{\begin{equation}}
\def\ee{\end{equation}}
\def\nn{\nonumber}
\def\bea{\begin{eqnarray}}
\def\eea{\end{eqnarray}}

\newtheorem{theorem}{Theorem}[section]
\newtheorem{lemma}[theorem]{Lemma}
\newtheorem{convention}[theorem]{Convention}
\newtheorem{remark}[theorem]{Remark}
\newtheorem{definition}[theorem]{Definition}

\DeclareMathSymbol{\varnothing}{\mathord}{AMSb}{"3F} 

\begin{document}

\title{Rotor interaction in the annulus billiard}

\author{P\'eter B\'alint} \author{Serge Troubetzkoy}

\address{Centre de Physique Th\'eorique and\\
Federation de Recherches des Unites de Mathematique de Marseille\\
CNRS Luminy, Case 907, F-13288 
Marseille Cedex 9, France\\ 
on leave from Alfr\'ed R\'enyi Institute of Mathematics 
of the H.A.S., H-1053, Re\'altanoda u. 13--15, Budapest, Hungary\\ 
and from Mathematical Institute, Technical University of Budapest, 
H-1111 Egry J\'ozsef u. 1., Budapest, Hungary}
\email{balint@cpt.univ-mrs.fr, bp@renyi.hu}
\urladdr{http://www.renyi.hu/{\lower.7ex\hbox{\~{}}}bp/}

\address{Centre de Physique Th\'eorique\\
Federation de Recherches des U\-ni\-tes de Mathematique de Marseille\\
Institut de math\'ematiques de Luminy and\\
Universit\'e de la M\'editerran\'ee\\ 
Luminy, Case 907, F-13288 Marseille Cedex 9, France}

\email{troubetz@iml.univ-mrs.fr}
\urladdr{http://iml.univ-mrs.fr/{\lower.7ex\hbox{\~{}}}troubetz/} \date{}
\subjclass{} 
\begin{abstract}
Introducing the rotor interaction in the integrable system of the 
annulus billiard produces a variety of dynamical phenomena, from
integrability to ergodicity.
\end{abstract}
\keywords{Billiards, rotor interaction, integrability, ergodicity, Anzai skew product.}
\maketitle

\pagestyle{myheadings}

\markboth{ROTOR INTERACTION IN THE ANNULUS BILLIARD}{P\'ETER B\'ALINT AND SERGE
TROUBETZKOY}

\section{Introduction}\label{sec1}
\setcounter{equation}{0}

The rotor model (see \eqref{evgen} below) has been introduced in 
\cite{MLL} and further studied from the physical point of view in \cite{LLM}.
The main motivation for replacing the usual elastic reflection in 
billiard models with such an interaction is that this way 
the a priori independent motion of point particles can be coupled. 
Despite of the physical interest in the subject, 
there is only a very little amount of 
mathematically rigorous discussion. This is related to the difficulty of the 
problem -- even the simplest possible {\it hyperbolic} candidates, planar 
Sinai billiards with circular scatterers and the related Lorentz gas models 
turned out to be rather complicated. In \cite{Bu} a detailed analysis has 
been given for a zero measure set of trajectories that possess special 
symmetries.

In the present article we investigate the rotor interaction in the context
of {\it elliptic} billiards. Consider the motion of point particles in 
the annulus billiard, i.e.~in the planar region bounded by two concentric 
circles. The one particle system is well known to be integrable (\cite{ann} p.~26),
while the dynamics of several point particles is the direct product of the 
individual regular systems. However, if we replace the inner (and possibly 
the outer) circular billiard scatterer with a rotating object, regularity 
may be destroyed. Three different models are investigated in this article. 
In Section~\ref{sec2} we consider the one particle system with the inner 
circular scatterer rotating and the outer one kept to be a billiard. For  
this system integrability persists. In Section~\ref{sec3} both scatterers 
are assumed to rotate and this setting destroys integrability.
For a positive measure set 
of  parameters velocity evolution is minimal (Theorem~\ref{thbase}) and the 
two dimensional system of our point particle is ergodic 
(Theorem~\ref{thskew}). Finally, in Section~\ref{sec4} the motion of {\it two} 
point particles is discussed in case the inner scatterer is rotating and the 
outer one is kept to be a billiard. Though this problem is 
significantly more complicated, 
it is possible to find an open set of parameters and, 
correspondingly, a positive measure set of initial conditions, that this 
system behaves as the one discussed in Section~\ref{sec3} for an arbitrary 
long finite time (Theorem~\ref{two} and Remark~\ref{twoskew}). 
From a physical point of view this 
result is probably the most interesting as it indicates, as claimed in 
the physics literature (see \cite{MLL}), the following phenomena: the rotor 
interaction switches on non-trivial coupling among the particles that 
would follow independent (in our case integrable) patterns if they were 
studied in the original billiard setting.

Throughout the paper we consider the motion of point particles in the annular 
region bounded by two concentric circles, the outer having radius one, the 
inner radius $R<1$. Interaction with the circular scatterers is either 
an elastic reflection (as in billiard models) or the scatterer is assumed
to rotate with some angular velocity $\omega$ and interact with the point 
particle according to the laws
\bea
\label{evgen}
v_n' &=& - v_n; \nn \\
v_t' &=& v_t -\frac{2\eta}{1+\eta}(v_t-r\omega);  \\
r\omega'&=& r\omega + \frac{2}{1+\eta} (v_t-r\omega).\nn 
\eea

Here the coordinates $v_n$ and $v_t$ are the normal and tangential components
(with respect to the circle of collision) 
of the particle velocity and the primes indicate post-collisional values. 
Furthermore $r$ is 
$1$ for the outer and $R$ for the inner circular scatterer, and 
$\eta=\frac{\Theta}{mr^2}$, where $\Theta$ is the moment of inertia for 
the scatterer in question and $m$ is the mass of the point particle. 

\begin{remark}
\label{notthesame}
Note that for the same velocity vector $\v$ the splitting into tangential 
and normal components is different for the inner and the outer scatterer. In 
Sections~\ref{sec2} and~\ref{sec4} it is only the inner scatterer that is 
assumed to rotate, and thus the components $v_n$ and $v_t$ in these sections 
always refer to the splitting at the inner scatterer. In Section~\ref{sec3} 
both scatterers rotate and it is more convenient to use different notation,
see the beginning of Section ~\ref{sec3.1}.
\end{remark}

\begin{convention}\label{conv}
The evolution of the coordinates in Formula \eqref{evgen} will be referred to
as the {\em base}.  
\end{convention}
In each section the dynamics of the base will be studied before turning
to the whole system, which is a skew product over the base, sometimes
referred to as fibers.

\section{One particle, one rotating scatterer}\label{sec2}
\setcounter{equation}{0}

This section is devoted to a model which is closest to the original integrable
annulus billiard, namely the motion of one point particle with only the inner
scatterer rotating.
For certain initial conditions the point particle never reaches the inner 
circular scatterer. In this invariant set motion is that of a point particle 
in a circle, a system well-known to be integrable. We consider the complement. 

At the inner scatterer we have
\bea
\label{ev1}
v_n' &=& - v_n; \nn \\
v_t' &=& v_t -\frac{2\eta}{1+\eta}(v_t-R\omega);  \\
R\omega'&=& R\omega + \frac{2}{1+\eta} (v_t-R\omega) \nn
\eea
which is just the inner version of Equation \eqref{evgen}, see Remark \ref{notthesame}.
\begin{lemma}
\label{lemma1}
(1) If a trajectory collides once with the inner scatterer, it collides with 
it 
infinitely often. More precisely, there is an alternating sequence of 
collisions, every first with the inner and every second with the outer 
scatterer.\\
(2) There are three integrals of motion.\\
(3) The velocity motion is periodic of period 2.\\ 
\end{lemma}

\begin{proof}
A trajectory which starts from the inner circle is reflected  elastically 
in the outer circle, thus it returns to the inner circle.
Iterating this argument proves 
statement (1).

Consider two consecutive reflections at the inner scatterer with one 
billiard reflection at the outer one in between. Let us denote by the triples
$v_n(1),v_t(1), \omega(1)$ and $v_n(2),v_t(2), \omega(2)$ the incoming 
velocities at the two consecutive inner reflections, and let us refer to the 
outgoing counterparts by primes. By the geometrical picture sketched above:
$$
v_n(2)=-v_n(1)'; \qquad v_t(2)=v_t(1)'
$$
and furthermore $\omega(1)'=\omega(2)$ as the inner scatterer does not 
interact while the particle collides outside. This, together with \eqref{ev1}, 
gives, on the one hand
$$
v_n(1)=v_n(2)=-v_n(1)'=-v_n(2)'
$$
which can be regarded as an integral of motion. On the other hand, to get   
$v_t(2)'$ and $\omega(2)'$ we need to iterate the second and the third 
formulas in \eqref{ev1} twice. This gives 
$$
\omega(2)'=\omega(1); \qquad v_t(2)'=v_t(1)
$$     
thus statement (3) is proved. 

Finally multiply the third equation in 
(\ref{ev1}) with $\eta$ and add it to the second to get
\be
\label{mom1}
v_t+\eta R \omega = v_t'+ \eta R \omega'=N.
\ee
Furthermore, subtracting the third line from the second in \eqref{ev1} we get
$$
v_t'-R\omega'= - (v_t-R\omega)
$$
which, together with (\ref{mom1}) results in
\be
\label{en1}
v_t^2+\eta R^2 \omega^2 = v_t'^2+  \eta R^2\omega'^2 = E.
\ee
Thus we have calculated the other two integrals.\end{proof}

\begin{remark}
Equation \eqref{mom1} is a multiple of the angular momentum, while 
\eqref{en1} is the ``tangential'' part of the full kinetic energy 
(the normal part of the kinetic energy, i.e.~the length of $v_n$ is 
in itself conserved as the third integral).

In the space of velocities $(v_t,\omega)$, that can be regarded as $\R^2$, 
fixing the values \eqref{mom1} and \eqref{en1} we get a line and a circle, 
respectively. These two curves intersect in (at most) two points 
$(v_t,\omega)$ and $(v_t',\omega')$. Velocity evolution is a period two cycle 
on these two points. 
\end{remark}

\begin{figure}[hbt]
\centering
\psfrag*{a}{\footnotesize$\hb$}
\psfrag*{b}{\footnotesize$\beta$}
\psfrag*{a'}{\footnotesize$\hb'$}
\psfrag*{b'}{\footnotesize$\beta'$}
\psfrag*{1}{\footnotesize$1$}
\psfrag*{R}{\footnotesize$R$}
\epsfxsize=10cm
\epsfysize=6cm
\epsfbox{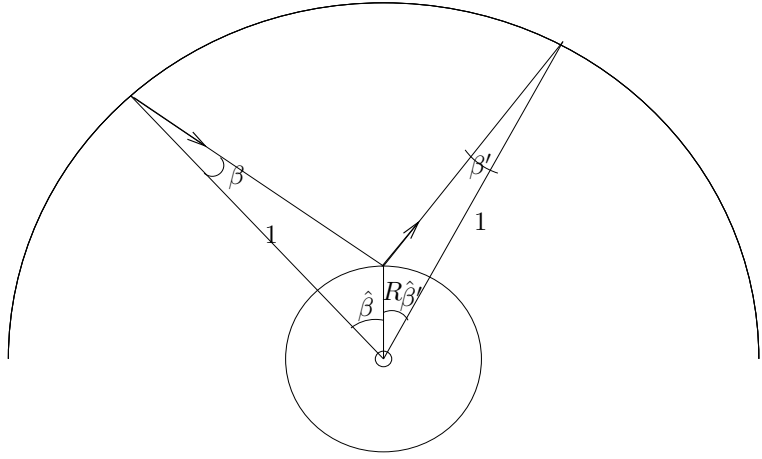}
\caption{$\hb={\rm arcsin}\left( \frac{\sin\beta}{R} \right) -\beta$}
\label{fig1}
\end{figure}

Turning to the evolution in the fibers, we essentially get the 
behavior observed in the integrable annulus billiard. Let us parameterize
the outer circle by the arclength parameter $\phi$ and denote the points of
consecutive impacts on a trajectory by $\phi_k, k\in\Z$. Then
$$
\phi_{k+1}=\phi_k+\hat{\beta}+\hat{\beta}'
$$
where
$$
\beta = {\rm arctan} \frac{v_t}{v_n}, \quad  \hat{\beta}={\rm arcsin}\left( \frac{\sin\beta}{R} \right) -\beta
$$
and the same for the primed quantities, see Figure~\ref{fig1}. 
Thus the angles $\hat{\beta},\hat{\beta}'\in 
[- {\rm arccos} (R),{\rm arccos} (R)]$ 
depend 
smoothly and strictly monotonically on the ratios $\frac{v_t}{v_n}$
and $\frac{v_t'}{v_n'}$, respectively.  
The consecutive $\phi_k$-s follow an orbit of 
a circle rotation, which is irrational apart from a countable collection
of one codimensional manifolds in the three dimensional parameter space of 
integrals. To see this keep the values of $E$ and $N$ (\eqref{en1} and 
\eqref{mom1}) fixed 
and increase the value of $v_n$ continuously, then $\hat{\beta}+
\hb'$ decreases 
continuously and is rational only for a countable set of $v_n$-s.

\section{One particle, two rotating scatterers}\label{sec3}
\setcounter{equation}{0}

\subsection{Coordinates, parameters and integrals of motion}\label{sec3.1}
In this section we will find that allowing both scatterers to rotate
results in a behavior different from the integrable models discussed
so far.

As long as velocities are concerned we need to handle an a priori 
four dimensional system: the coordinates are, on the one hand,  
$\omega_1,\omega_2$ for the angular velocities of the inner and outer 
scatterers, respectively, and, on the other hand,  
$\v$, the two dimensional velocity vector of 
the point particle.
Recall Remark~\ref{notthesame}: in contrast to the previous section, we 
will use notations $v_n$ and $v_t$ for the splitting of the velocity
$\v$ into normal and tangential components at the {\it outer} scatterer. 
The sign of $v_t$ has a physical meaning, it describes the direction of
the ``circular'' component of the motion of the point particle.
On the contrary, the sign of $v_n$ 
is in a certain sense irrelevant, the positivity/negativity of $v_n$ merely
tells
us if we are just before/after a collision with the outer circle, thus we
fix $v_n$ to be always positive.
The symbols $\bv_n,\bv_t$ are used for the splitting of $\v$ into normal 
and tangential components at the inner scatterer. By simple geometrical 
observation they can be defined for vectors 
$\v$ which satisfy $\frac{v_t^2}{v_n^2}\le \frac{R^2}{1-R^2}$ as
$\bv_t=\frac{1}{R} v_t$ and $\bv_n=\sqrt{v_n^2-\frac{1-R^2}{R^2}v_t^2}$. 
Otherwise the trajectory starting from the outer scatterer reaches it again
without interacting with the inner one. 

We use the notation
$\omega_i$ for the angular velocity and $\Theta_i$ for the moment of inertia 
with $i=1,2$ in the case of the inner and the outer scatterer, respectively.
The rescaled moments of inertia are $\eta_1=\frac{\Theta_1}{mR^2}$ and 
$\eta_2=\frac{\Theta_2}{m}$. In these coordinates,  the interaction at the 
outer scatterer is given by:
\bea
\label{out}
v_t' &=& v_t -\frac{2\eta_2}{1+\eta_2}(v_t-\omega_2); \nn \\
\omega_2'&=& \omega_2 + \frac{2}{1+\eta_2} (v_t-\omega_2). 
\eea
By the convention on the sign of $v_n$ made above, the normal velocity is 
preserved and straightforward calculation gives two further preserved  
quantities:
\be
\label{intout}
v_t+\eta_2 \omega_2=v_t'+\eta_2 \omega_2'; \qquad 
v_t^2+\eta_2\omega_2^2 = v_t'^2+\eta_2\omega_2'^2. 
\ee
As to interaction with the inner scatterer
\bea
\bv_t' &=& \bv_t -\frac{2\eta_1}{1+\eta_1}(\bv_t-R\omega_1); \nn \\
R\omega_1'&=& R\omega_1 + \frac{2}{1+\eta_1} (\bv_t-R\omega_1).\nn 
\eea
Note, however, that $v_t=R\bv_t$ and thus these equations are equivalent to
\bea
\label{inn}
v_t' &=& v_t -\frac{2\eta_1}{1+\eta_1}(v_t-R^2\omega_1); \nn \\
R^2\omega_1'&=& R^2\omega_1 + \frac{2}{1+\eta_1} (v_t-R^2\omega_1). 
\eea
Two preserved quantities for the inner collisions are
\be
\label{intinn}
v_t+\eta_1 R^2 \omega_1=v_t'+\eta_1 R^2 \omega_1'; \qquad 
v_t^2+\eta_1 R^4 \omega_1^2 = v_t'^2+\eta_1 R^4 \omega_1'^2. 
\ee
As $\omega_1$ ($\omega_2$) does not change at outer (inner) collisions, 
\eqref{intout} and \eqref{intinn}
imply the presence of two integrals of motion, the angular momentum
\be
\label{mom}
v_t+ \eta_1 R^2 \omega_1 + \eta_2 \omega_2 ={\rm const.} = N
\ee
and the tangential energy
\be
\label{en}
v_t^2+ \eta_1 R^4 \omega_1^2 + \eta_2 \omega_2^2={\rm const.} = E (>0).
\ee
The term ``tangential'' refers to the fact that the 
full kinetic energy of the system is a third integral of motion, always 
greater than the tangential one. More precisely, dividing the full energy 
$\frac{1}{2} (m\v^2+\Theta_1\omega_1^2+\Theta_2\omega_2^2)$ by
$\frac{1}{2}m$ we get
\be
\label{full}
v_n^2 + v_t^2 + \eta_1 R^2 \omega_1^2  + \eta_2 \omega_2^2 ={\rm const.} 
= F (>E).
\ee
Note that \eqref{intout} (\eqref{intinn}) and \eqref{full} together 
imply that $v_n$ ($\bv_n$) is preserved during inner (outer) collisions.  

\begin{remark}\label{novzero}
(1) One of our aims below is to understand how the dynamical behavior depends 
on the several parameters of the system. The constants 
$\eta_1$ and $\eta_2$ will be 
referred to as physical parameters in contrast to the integrals of motion 
$N, E$ and $F$.

(2) To avoid complications we assume $\v \ne \vec{0}$. Note that this assumption is
reasonable in the following sense: $\v=\vec{0}$ is only possible for a 
codimension one submanifold of integrals of motion. To see this observe that 
if we had $\v= \vec{0}$, \eqref{mom} and \eqref{en} would fix the values of the 
other two coordinates $\omega_1$ and $\omega_2$ (up to two possibilities) 
and thus would create, via \eqref{full}, a relation between the three integrals
of motion. 
\end{remark}

\subsection{Evolution of velocities}
\label{secvel}

For convenience we rescale our velocities as 
\be
\label{rescale}
x:=\sqrt{\eta_1}R^2 \omega_1,\ y:=\sqrt{\eta_2}\omega_2, \ z:=v_t, \ w:= v_n.
\ee
In these rescaled quantities our integrals of motion are
\be
\label{intnew}
x^2+y^2+z^2=E; \qquad  \sqrt{\eta_1} x + \sqrt{\eta_2} y +z = N
\ee
and
\be
\label{fullnew}
\frac{x^2}{R^2} + y^2 + z^2 + w^2 = F.
\ee
By \eqref{out} and \eqref{inn} we can first consider how the velocity 
coordinates $x,y,z$ evolve at collisions. With the convention on its 
sign the value of $w(=v_n)$ is determined by \eqref{fullnew} 
for each fixed value of the other three coordinates. 
In the rest of the section, unless 
otherwise stated, velocity evolution will always refer to the dynamics in 
the triple of $x,y,z$. 

Taking into account that $x$ ($y$) 
does not change at outer (inner) collisions, 
in the rescaled coordinates the transformations \eqref{out} and \eqref{inn}
are reflections across the planes
\be
\label{planes}
\sqrt{\eta_2}z-y=0 \quad {\rm and} \quad \sqrt{\eta_1}z-x=0,
\ee
respectively. Furthermore, by
\eqref{intnew} motion is restricted to the intersection of a sphere and a 
plane, i.e.~to a circle which will be referred to as the velocity circle.
We will see that the base, see Convention \ref{conv}, will consist of two copies of this circle.
Note that these circles do not describe the physical position of the particle,
on the other hand they give the state of the velocities.
On the velocity circle the reflections across the planes \eqref{planes} reduce
to reflections across lines. For the outer collisions we have reflections 
across
\be
\label{line2}
l_2(t)= \left( \, -\frac{1+\eta_2}{\sqrt{\eta_1}}t\, ,\, \sqrt{\eta_2}t\, ,
\, N+t\, \right) 
\ee  
while for the inner ones across
\be
l_1(t)= \left( \, \sqrt{\eta_1}t\, , \, -\frac{1+\eta_1}{\sqrt{\eta_2}}t\, ,
\, N+t\, \right).\nn 
\ee  

Let us denote the angle between these two lines by $\frac{\gamma}{2}$. 
Depending on it we distinguish rational, irrational (and later on even 
Diophantine) cases.

\begin{lemma}
\label{baseangle}
Rationality is independent of the integrals of motion. On the other hand, 
there exists a countable set of codimension 1 submanifolds of physical
parameters which lead to the rational case, otherwise (and thus 
for a full Lebesgue measure set of parameters) $\gamma$ is irrational.  
\end{lemma}
\begin{proof}
Straightforward 
calculation gives
\be
\label{gamma}
\cos \frac{\gamma}{2} = \left( (1+{\eta_1}^{-1})(1+{\eta_2}^{-1}) 
\right)^{-\frac{1}{2}}.
\ee
Thus the rationality of $\gamma$ depends only on the physical parameters. 
Consider the numbers $\cos(\frac{p}{q}\frac{\pi}{2})$ for positive integers
$p<q$. Note that elements of this countable set $A$ are all algebraic numbers. 
The base is rational iff the left hand side of \eqref{gamma} is equal to 
some number in $A$. This happens for a countable collection of one 
codimensional submanifolds of physical parameters.   
\end{proof}

\begin{convention}\label{01} To agree with different standard
conventions, in
the geometric framework we will think 
of $\S^1$ as $[0,2\pi)$, while in the number theoretic framework
we think of $\S^1$ as $[0,1)$. In both cases we scale all $\S^1$ 
valued quantities ($s,\phi,\gamma,\alpha(s)$) appropriately. 
This should cause no confusion.
\end{convention}

To give a description of the base dynamics we introduce the arclength 
parameter $s$  as a coordinate on the velocity circle. For definiteness 
$s=0$ corresponds to one of the intersection points of the circle with 
the line \eqref{line2}. In this coordinate the outer and inner bounces, 
i.e.~the transformations \eqref{out} and \eqref{inn} are
\bea
s &\to& s_1 = -s \  ({\rm mod} \, 1) {\rm \  and}\label{outs}\\ 
s &\to& s_2 = -s-\gamma \ ({\rm mod} \, 1) \label{inns},
\eea  
respectively. 
Take two copies, $O$ and $I$, of the velocity circle describing the 
outgoing velocities just after the
bounces at the inner and outer circles, respectively. The phase space 
of the velocity motion, or the base, 
is
$$
\hat{M}=\S^1 \times \{ 1,2 \}= I \cup O; \qquad \hat{M}\ni \bs=(s,i) 
\ {\rm with} \  s\in\S^1 \ {\rm and} \ i=1,2.
$$
The velocity evolution will be denoted by $T: \hat{M}\to \hat{M}$.
As an inner bounce is always followed by an 
outer one, for $\bs \in I$, i.e.~for $(s,1)$ the image $T\bs$ is always in $O$,
more precisely $T(s,1)=(s_1,2)$ (recall \eqref{outs}).
An outer bounce, however, may be followed by either an inner or an outer 
bounce, depending on the angle the point particle velocity makes with the
normal vector of the outer scatterer.
\begin{lemma}
\label{alter}
(1) There is a (possibly empty) open set $U\subset \S^1$, such that for 
$\bs\in U$ (more precisely for $\bs=(s,2)$ with $s\in U$), we have 
$T\bs=(s_1,2)$ (the image is the point defined by \eqref{outs} in $O$), 
and for $\bs\in O\setminus U$ (more precisely for $\bs=(s,2)$ with 
$s\in \S^1\setminus U$) we have $T\bs=(s_2,1)$ 
(the image is the point defined by \eqref{inns} in $I$).\\
(2) The set $U\subset \S^1$ is invariant with respect to the reflection 
\eqref{inns}.\\
(3) There is an open set of integrals of motion for which $U$ is empty.
\end{lemma}

\begin{remark}
\label{alt}
We will refer to empty $U$ as the alternating case, non-empty $U$ as 
the non-alternating case. Both of them occur for open sets of 
integrals of motion, independently of the value of physical parameters. 
\end{remark} 

\begin{proofof}{Lemma~\ref{alter}}

(1) By Formulas \eqref{intnew} and \eqref{fullnew}
the arclength $s$ on the velocity circle determines $z=v_t$ and 
$w=v_n\ge 0$ continuously.
The geometric condition for $(s,2)=\bs\in O$ to reach the outer scatterer 
without 
touching the inner one (and thus interacting with it) is 
$\frac{v_t^2}{v_n^2}>\frac{R^2}{1-R^2}$. Note that at tangential collisions 
with the inner scatterer $v_t$ is modified by \eqref{inn}, thus this 
corresponds to landing on $I$. 

We define
\be
\label{U}
U:= \left\{ \ s\in \S^1 \ | \ z^2(s)>\frac{R^2}{1-R^2} w^2(s)\  \right\} .
\ee
If $\bs\in O$ with $s\in U$ (i.e.~in case 
the next bounce is outside) then \eqref{outs} applies, and thus 
$T\bs=(s_1,2)$. 
Otherwise (i.e.~if $\bs\in O\setminus (U\times\{ 2 \})$)
the next bounce is inside, thus the image is the point defined by
\eqref{inns} in $I$, or $T\bs=(s_2,1)$. 
The set
$U\subset \S^1$ is the preimage of an open set by a continuous map, thus it 
is open. 

(2) To see that $U$ is invariant under \eqref{inns} note first that $y$ does 
not change during an inner collision. This, together with the presence of the 
integrals \eqref{intnew} and \eqref{fullnew} implies that $w^2- 
\frac{1-R^2}{R^2}z^2$ is invariant for inner bounces, i.e.~for the 
reflections \eqref{inns}. The statement follows as it is exactly the sign of 
this quantity that distinguishes the points of $U$.

(3) The set $U$ is empty if 
\be
\label{emptyu}
z^2(s) <  \frac{R^2}{1-R^2} w^2(s)
\ee 
for all $s$. 
Consider the integrals of motion defined by \eqref{intnew} and 
\eqref{fullnew} and denote $\lambda=\frac{F}{E}>1$. 
By \eqref{intnew} $x^2$ and $z^2$ are both less than or equal to 
$E$ for all $s$. This, together with \eqref{fullnew} gives
\bea
w^2 &\ge&  F-y^2-z^2-x^2 - x^2 ( R^{-2}-1 ) \, \ge \, F - E - (R^{-2} -1)E 
\, = \,\nn \\ 
&=&  F- ER^{-2} \, = \, E (\lambda-R^{-2}) \, \ge \, 
(\lambda - R^{-2}) z^2 \nn
\eea  
for all $s$. This implies $U$ is empty if $\lambda -R^{-2} >
\frac{1-R^2}{R^2}$, i.e.~if  $\lambda > (\frac{2}{R^2}-1)$, a condition 
satisfied for a set 
of parameters that has nonempty interior.  
\end{proofof}

\begin{remark}
\label{nouin}
(1) Actually $U$ is empty iff we replace the strict inequality in 
Formula~\eqref{emptyu}
with a nonstrict one.  The analysis of the equality case is substantially more complicated,
and is a co-dimension one phenomenon in the set of parameters, thus we do not analyze it.

(2) Note that the points $s\in U$ correspond to velocity configurations that 
cannot be realized at the inner circle. Thus for nonempty $U$ the phase 
space is $M\subset \hat{M}$:
\be
M=I \cup O =\left(  (\S^1\setminus U)\times\{1 \} \right) \cup 
\left( \S^1\times \{ 2 \} \right)\nn 
\ee
and the dynamics acts as
\bea
T \big{(} O\setminus(U\times \{2 \}) \big{)} &=& 
(\S^1 \setminus U)\times \{ 1\} = I; \nn \\
T(U\times \{ 2\}) \cup TI &=& \S^1\times \{ 2 \} = O;\nn 
\eea
with
\bea
T(s,2) &=& (s_2,1) \quad  {\rm for} \quad  s\in \S^1\setminus U; \nn \\  
\label{realdin}
T(s,2) &=& (s_1,2) \quad {\rm for} \quad  s\in U;  \ {\rm and}\\ \nn
T(s,1)&=& (s_1,2) \quad {\rm for} \quad  s\in \S^1\setminus U;
\eea
which is well-defined and invertible 
by statement (2) of Lemma~\ref{alter}. As $U$ has no 
inside copy sometimes we refer to $U\times \{ 2\}$ as $U$ for brevity, 
this should cause no confusion.
\end{remark}

The main result on the base dynamics is the following theorem.
\begin{theorem}
\label{thbase}
In the alternating irrational case the dynamics of the base is minimal, while 
in all the other cases velocity motion is always periodic. 
Moreover, all possible periods are even with length bounded above by a 
constant  (which depends on the parameters).  
\end{theorem}

\begin{proof}
Consider the alternating case first. For $\bs\in O$ all odd iterates are in 
$I$ and all even iterates are in $O$ (and for $\bs \in I$ the other way 
round). Thus it is enough to consider the first return map onto $O$ (which 
is exactly $T^2$). By \eqref{outs} and \eqref{inns} this is a rotation of 
$\S^1$ by $\gamma$, which is periodic/minimal exactly in the 
rational/irrational cases. (The odd iterates make the reflected image of 
the same rotation orbit in $I$.)

Turning to the non-alternating case, i.e.~to $U\ne \emptyset$, first note 
that $U\cap T^{-1}U$ is an invariant set. To see this consider 
$(s,2)=\bs\in U\cap T^{-1}U$. Then both $T\bs$ and $T^2\bs$ are in $O$. 
More precisely
$$
T(s,2)=(s_1,2) \quad {\rm and} \quad T^2(s,2)=T(s_1,2)=(s,2),
$$  
thus the set is invariant, furthermore, all points of it are periodic 
of period $2$.

Now take $(s,2)=\bs\in U\setminus T^{-1}U$. Then $T\bs$ is in $O\setminus U$,
$T^2\bs \in I$ and even $T^{-1}\bs\in I$, as the point is assumed not to 
belong to the above described invariant set. Moreover,  
by \eqref{outs} and \eqref{realdin}
\be
\label{turnback1}
T\bs=T(s,2)=(s_1,2) \quad {\rm and} \quad T^{-1}\bs= T^{-1}(s,2)=(s_1,1),
\ee
thus the trajectory essentially ``turns back''. The time reflection symmetry
of \eqref{turnback1} is preserved until the next $U$-hitting time. To see this
let $k>1$ be the smallest integer (possibly $\infty$) such that 
$T^k\bs\in U$. Up to time $k$ the orbit is alternating, i.e.~by induction
\bea
T^{2l-1}\bs=(s',2), \qquad T^{2l}\bs=(s'_2,1), \nn \\
T^{-(2l-1)}\bs=(s',1),
\qquad T^{-2l}\bs=(s'_2,2) \nn 
\eea
for some $s'\in \S^1\setminus U$ whenever $0<2l<k$. In particular 
$T^{(k-1)}\bs=(s'',1)$ and $T^{-(k-1)}\bs=(s'',2)$ for some
$s''\in T^{-1}U\setminus U$ (note that $k$ is always odd). This implies by 
\eqref{realdin} that the preimage of $T^{-(k-1)}\bs$ is not in $I$, 
alternation ceases and
$$
T^{-k}\bs=(s_1'',2)\in U \subset O.
$$
Moreover 
$$
T^k\bs=(s_1'',2)=T^{-k}\bs,
$$
thus the point is periodic with period $2k$. This way we have shown that any 
orbit which enters $U$ twice is periodic.

Consider any $\bs=(s,2)\in O$ and the subset 
of $\S^1$ generated from $s$ by the 
reflections \eqref{outs} and \eqref{inns}. This subset of $\S^1$ consists
of two rotation orbits: the first by $\gamma$ and the second by $-\gamma$,
 which is the reflected image  of the first one by \eqref{inns}.  

In the non-alternating rational case there are two possibilities. On 
the one hand, if both of the above mentioned two rotation orbits avoid $U$,
dynamics is that of the alternating case: $\bs$ is periodic of length that 
depends only on $\gamma$ (i.e.~on the physical parameters). On the other 
hand, by the invariance of $U$ with respect to \eqref{inns} (cf.~statement (2) 
of Lemma~\ref{alter}) if one of the rotation orbits hit $U$, the other orbit 
-- the reflected image -- should hit it as well. Thus if this occurs the 
trajectory of $\bs$ enters $U\times \{ 2 \}$ twice and becomes 
periodic with even length. Since the longest possible period
occurs when we avoid $U$,  the length of all possible periods 
is bounded above by a 
constant that depends only on $\gamma$.

In the non-alternating irrational case, as both rotation orbits are dense and 
thus cannot avoid $U$, the second of the two possibilities sketched above
occurs for all $\bs$. To 
complete the proof of the theorem we only need to show that the length of
possible periods is bounded above. By minimality and compactness there exists 
a positive integer $N$ such that the first $N$ rotation-preimages of 
$U$ cover $\S^1$. This implies that the length of any period is less than or 
equal to $4N$. Moreover, $N$ depends only on $\gamma$ and $U$, i.e.~on the 
physical parameters and the integrals of motion.
\end{proof}

\subsection{Analysis of the skew product}
\label{secskew}

In this subsection we extend our investigation to the fibers, i.e.~to the 
question how the points of impact for the consecutive outer bounces follow 
each other on the outer scatterer. Let us parameterize the outer circle by 
the arclength parameter $\phi$ (the starting point is irrelevant and 
may be chosen arbitrarily).  We consider the skew product
\be
\label{skew}
F:\S^1 \times \S^1 \to \S^1 \times \S^1; \qquad F(s,\phi)=
(T_O s, \phi +\alpha(s)).
\ee
Here $T_O$ is the (projection onto $\S^1$ of the) first return map onto $O$ 
for $T$ from Section~\ref{secvel}. Note that $T_O$ is either the first or 
the second iterate of $T$. More precisely, if $s\in U$, we have $T_O=T$ and 
$\alpha(s)=\pi -2\beta(s)$ where $\tan \beta(s)=\frac{v_t(s)}{v_n(s)}=
\frac{z(s)}{w(s)}$. On the other hand, if $s\in \S^1\setminus U$, we have 
$T_O=T^2$ and $\alpha(s)$ can be calculated using the geometry of 
Figure~\ref{fig1} (cf. the end of Section~\ref{sec2}) as:
\bea
\alpha(s) &=& \hat{\beta}(s)+\hat{\beta}(s_2) \qquad {\rm where} \nn \\ 
\label{alphas}
\sin(\beta(s)+\hat{\beta}(s)) &=& \frac{\sin\beta(s)}{R}, \\ \nn
\tan \beta(s) &=& \frac{v_t(s)}{v_n(s)}=\frac{z(s)}{w(s)}
\eea 
and $s_2$ is defined by \eqref{inns}.
Note that, when restricted to $U$ or $O\setminus U$, both
$\alpha$ and $T_O$ are continuous.

Recall from Theorem~\ref{thbase} that apart from the alternating irrational 
case $T_O$ is periodic for any $s$ with some period $N=N(s)$. In these 
periodic situations
$$
F^N(s,\phi)=(s,\phi+\alpha_N(s)); \quad \alpha_N(s)= \alpha(s)+\alpha(T_O s)
+\cdots+\alpha({T_O}^{N-1}s).
$$
Thus a suitable iterate of the dynamics in the fibers is a rotation, 
the angle of which depends only on the initial velocity $s$, and not on the
initial position $\phi$. The iterates $F^{k_0+kN}s$, with $0<k_0<N$ fixed and 
$k\in\Z$, are the shifted versions of the same rotation orbit. Moreover, by 
piecewise continuity of $\alpha$ and $T_O$, and by the uniform upper bound on 
$N(s)$, the rotation angle $\alpha_N(s)$ depends on $s$ in a piecewise 
continuous manner as well (at discontinuity points, however, even $N(s)$ may 
jump).

\begin{remark}
Note that in the piecewise continuity statements we have finitely many pieces
since the set $U$,being defined by algebraic equations and inequalities, 
has finitely many connected components.
\end{remark}

From this point on we restrict to the alternating irrational case.
Thus $T_O$ is an irrational rotation, and we 
have
\be
\label{skewf}
F (s,\phi) = (s+\gamma,\phi+\alpha(s))
\ee
with $\gamma$ irrational and $\alpha:\S^1\to \S^1$ continuous. 
Such systems are 
called Anzai skew products after their first investigator (\cite{A}). 
There is an extensive literature on the subject, see \cite{F},\cite{BM}, 
\cite{I1}, \cite{I2} and references therein. 
Some simple issues are even treated in 
general monographs on dynamical systems like \cite{KH}.

\begin{figure}[hbt]
\centering
\psfrag*{db}{\footnotesize$d(\beta)$}
\psfrag*{b}{\footnotesize$\beta$}
\psfrag*{sb}{\footnotesize$\sin\beta$}
\psfrag*{1}{\footnotesize$1$}
\psfrag*{R}{\footnotesize$R$}
\epsfxsize=8cm
\epsfysize=8cm
\epsfbox{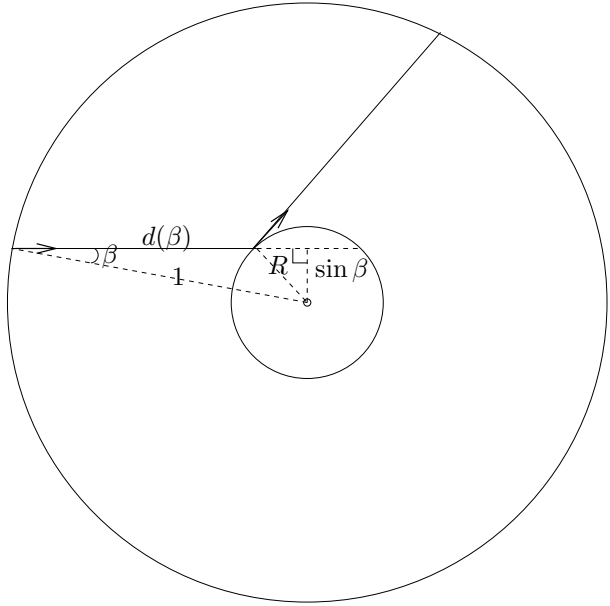}
\caption{$d(\beta)=\cos\beta -\sqrt{R^2-\sin^2\beta}$}
\label{fig2}
\end{figure}

In addition to the skew product map \eqref{skewf} 
we are interested in the dynamical 
behavior in real time, i.e.~in a suspension flow over our skew product. The 
ceiling function for this flow, i.e.~the physical (real-valued)
first return time to the outer scatterer depends only on the initial 
velocity, thus it may be denoted as $\tau(s)$ with 
$\tau:\S^1 \to \R^+$. We may calculate $\tau(s)$ by 
Figure~\ref{fig2} as 
\bea
\tau(s)&=&\frac{d(s)}{|\v(s)|} + \frac{d(s_2)}{|\v(s_2)|} \qquad {\rm with} \nn \\
\label{Ts}
|\v(s)| &=& \sqrt{v_t(s)^2 + v_n(s)^2} = \sqrt{z(s)^2+w(s)^2} 
\quad  {\rm and} \\  \nn 
d(\beta(s)) &=& d(s) = \cos \beta(s) -\sqrt{R^2-\sin^2\beta(s)}, \quad  {\rm  where} \\ \nn
\tan \beta(s) &=& \frac{v_t(s)}{v_n(s)}=\frac{z(s)}{w(s)}
\eea 
and $s_2$ is defined by \eqref{inns}. 

We recall some basic concepts and properties related to 
irrational numbers $\gamma$ and continuous 
functions $\alpha:\S^1\to\S^1$ (on details see 
\cite{KH} and \cite{nu}).

\begin{definition}
\label{dio}
Given $r>1$ an irrational number $\gamma\in[0,1]$ is 
{\bf Diophantine} of type $r$ 
if there exists $c>0$ such that for any $p,q\in \Z$, $q\ne0$
\be
|q\gamma-p|>c q^{-r}. \nn 
\ee
A number which satisfies the above condition for some $r>1$ is called 
Diophantine. Liouville numbers $\gamma$ are those irrationals 
that do not satisfy the 
Diophantine condition for any $r$. It is not hard to see that the set of
Diophantine numbers of type $r$ is, as a subset of $[0,1]$, on the one hand, 
of the first category, 
and, on the other hand, of full Lebesgue measure, with, even more, 
its complement 
having Hausdorff dimension $r^{-1}$. Consequently, the set of Diophantine 
numbers is of full measure and first category, while its complement, 
the union of Liouville numbers and rational ones, has 
zero Hausdorff dimension.   
\end{definition} 

\begin{definition}
\label{degree}
For any continuous map $\alpha:\S^1\to \S^1$ there exist countably many 
continuous functions $A:\R\to \R$ such that $A(s+1)-A(s)$ is an integer 
independent of $s$ and 
$\alpha(s)$ can be regarded as the restriction of $A(s)$ to $[0,1)$. 
The integer $A(s+1)-A(s)$ 
is independent of the lift and is called  {\bf the degree} of 
$\alpha(s)$. The map $\alpha(s)$ can 
be regarded as a real-valued continuous function defined on $\S^1$ 
iff its degree is zero.  
\end{definition}

\begin{definition}
\label{cohom}
Two measurable maps $\alpha:\S^1\to \S^1$ and $\alpha':\S^1\to \S^1$ are 
said to be {\bf cohomologous} (for the rotation by $\gamma$) iff there exists
$V:\S^1\to\S^1$ measurable such that
$$
\alpha(s)=\alpha'(s)+V(s+\gamma)-V(s)
$$ 
for almost all $s\in\S^1$. Maps cohomologous to $0$ are called coboundaries.
\end{definition}

The following statements date back to \cite{A}.
\begin{lemma}
\label{anzai}
(1) The skew product \eqref{skewf} is ergodic (with respect to Lebesgue 
measure) iff there is no $p\in\Z\setminus \{0\}$ such that $p\alpha(s)$ is a 
coboundary.\\
(2) The skew product \eqref{skewf} has pure point spectrum iff
$\alpha(s)$ is not cohomologous to a constant.\\ 
(3) A constant function $\alpha(s)=\alpha$ is a coboundary iff $\alpha$ and 
$\gamma$ are rationally dependent.\\
(4) By (1-3) for a function $\alpha(s)$ cohomologous to a constant $\alpha$ 
the skew product is either not ergodic or ergodic with pure point spectrum, 
depending on whether $\alpha$ and $\gamma$ are rationally dependent or not. 
\end{lemma}

The following result is a consequence of Theorem 1 from \cite{BM} combined
with Lemma \ref{anzai}.

\begin{lemma}
\label{nomixdio}
Consider an Anzai skew product \eqref{skewf} with $\gamma$ Diophantine and 
$\alpha:\S^1\to\S^1$ of zero degree and $C^\infty$. Then there exists a number
$\alpha\in [0,1)$ such that the map $F$ is smoothly (i.e.~$C^{\infty}$) 
conjugate to 
\be
\label{doublerot}
F_0 (\S^1\times\S^1)\to(\S^1\times\S^1), \qquad 
F_0(s,\phi)=(s+\gamma,\phi+\alpha).
\ee
\end{lemma} 

We prepare for our main result related to the full system, 
Theorem \ref{thskew},  with the proof of 
the following Lemma.

\begin{lemma}
\label{cinfty}
(1) The map $\alpha(s)$ is of zero degree, thus can be regarded as 
$\alpha:\S^1 \to \R$. \\
(2) Both $\alpha(s)$ and $\tau(s)$ are $C^{\infty}.$
\end{lemma}

\begin{proof}
(1) Note that any continuous map $\S^1\to\S^1$ that has nonzero degree must 
be onto. On the other hand, as 
we are in the alternating case we have $\cos\hat{\beta}(s)>R$ for any 
$s\in \S^1$ (see Figure~\ref{fig1}). 
This means $|\alpha(s)|\le 2{\rm arccos} R$ which in particular implies 
$\alpha(s)$ cannot be onto, thus it is of zero degree. 

(2) The maps $\alpha(s)$ and $\tau(s)$ are defined in Formulas~\eqref{alphas} 
and \eqref{Ts} respectively (see also \eqref{inns}). To establish the desired
smoothness properties we need to show that four functions, namely   
(i) $\beta(s)$ (ii) $\v(s)$ (iii) $\hat{\beta}(\beta)$
(iv) $d(\beta)$ are all smooth.  Moreover we need to show that
(v) $|\v(s)|$ is bounded away from zero. 

First note $x,y,z$ all depend on $s$ smoothly, while $w$ a priori only 
continuously. To see this note that $x(s),y(s),z(s)$ are the Cartesian
coordinates of a circle (of the velocity circle, cf. Section~\ref{secvel})
thus depend smoothly on the arclength parameter. By \eqref{fullnew} $w(s)$ 
is the square root of a smooth non-negative function, thus it is continuous, 
and smooth if bounded away from zero, a fact we establish below.

By  Remark~\ref{novzero} (2) and compactness $|\v(s)|$ is bounded
away from zero, thus we have statement (v). 
This by alternating (cf.~\eqref{emptyu}) gives positive 
lower bound and, by the above reasoning, 
smoothness for $w(s)$ which means (cf.~\eqref{rescale}) we have shown
(ii). 
Moreover (i) follows as we have demonstrated that $\tan \beta(s)$ 
is smooth.

We need to prove (iii) and (iv). As we are in the alternating case,
see the first part of Remark \ref{nouin},
we have Formula \eqref{emptyu} and thus by compactness
\be
R^2-\sin^2\beta(s)>K_1>0; \qquad {\rm and} \qquad 
\frac{\sin^2\beta}{R^2}>K_2>1.\nn
\ee
This means that the quantity under the square root in the 
formula for $d(\beta)$ (cf.
\eqref{Ts}) is bounded away from zero. Similarly, in the formula
for  $\hat{\beta}(\beta)$ (see \eqref{alphas}), the argument of
arc sin is bounded away from $1$ and $-1$. 
Thus we have (iii) and (iv).
\end{proof}

\begin{theorem}
\label{thskew}
Let us fix $\gamma$ Diophantine. \\
(1) The map \eqref{skewf} is smoothly conjugate to \eqref{doublerot} with
$\alpha=\int_{\S^1}\alpha(s)ds$.\\
(2) There is a positive measure set of integrals of motion for which 
\eqref{skewf} is uniquely ergodic. \\
(3) The system in ``physical'' time, i.e.~the suspension flow over \eqref{skewf}
with ceiling function \eqref{Ts} is never ergodic.
\end{theorem}

\begin{proof}
(1) Lemmas~\ref{cinfty} and~\ref{nomixdio} imply that \eqref{skewf} is smoothly
conjugate to \eqref{doublerot}. 

(2) For $\alpha$ and $\gamma$ rationally independent the system 
\eqref{doublerot} is well-known to be uniquely ergodic,
see \cite{KH}. What we are going to demonstrate is that this
happens for a positive measure set of parameters.  To see this, first note
that the velocity circle is determined by the two parameters $N$ and $E$
(see Equation \eqref{intnew} and the discussion at the beginning of Section \ref{secvel}).
Furthermore, it is possible to choose these two parameters from a positive measure
set in $\R^2$ to ensure that $x(s) = v_t(s) > 0$ for all $s \in \S^1$.  
We keep the physical parameters (and thus $\gamma$) together with two integrals of motion, with 
$N$ and $E$ from the above set fixed and vary the third one, 
$F$ from \eqref{fullnew}. What is shown below is that in this setting 
$\alpha$ depends on $F$ in a continuous and strictly monotonic manner. This 
implies that $\alpha$ and $\gamma$ can be rationally dependent only for 
countably many values of $F$,
and thus by our choice of $N$ and $E$ for a positive measure set of integrals.

In the rest of the proof, all functions of $s$ 
(e.g.~$z(s),\beta(s),\alpha(s)$) refer to the original 
values of the integrals, while the primed ones 
(e.g.~$z'(s),\beta'(s),\alpha'(s)$) indicate maps for $F$ slightly 
increased. 
Since  Lebesgue measure on the circle is invariant with respect to the
transformation $s \to s_2$, we have
$$
\alpha=\int_{\S^1}(\hat{\beta}(s)+\hat{\beta}(s_2))ds=
2\int_{\S^1}\hat{\beta}(s)ds.
$$

Note that
$x(s)=x'(s)$ and $y(s)$, $z(s)$ similarly, however as $F$ is increased
we get $w'(s)>w(s)$ for all $s \in \S^1$. 
This implies $\beta'(s) < \beta(s)$ (recall our convention on
$x(s) > 0$) which yields, by a simple geometrical calculation,
$\hat{\beta}'(s) < \hat{\beta}(s)$ for all $s \in \S^1$.
Thus $\alpha$ depends on $F$ strictly monotonically.

(3) A  suspension flow above an ergodic map is ergodic iff the 
ceiling function is not cohomologous to a constant, see \cite{KH}. As
$\tau(s)$ is $C^{\infty}$ by 
Lemma~\ref{cinfty} and $\gamma$ is chosen to be Diophantine, we 
may apply Lemma~\ref{nomixdio} (note also Definition \ref{degree}) to see that $\tau(s)$ is 
smoothly 
cohomologous to its average, a constant, thus the suspension flow is not 
ergodic.
\end{proof}

\begin{remark}
\label{weakmix}
Note that an Anzai skew product \eqref{skew} is never mixing as it has a 
factor which is a rotation by $\gamma$ (with 
the semi-conjugacy being projection to the first $\S^1$-coordinate). On the 
other hand, in case the spectrum is not pure point only 
``horizontal'' sets $B'=B\times\S^1$ do not mix weakly, see Lemma \ref{anzai}.
To be
more precise, if the spectrum is not pure point, restricting the unitary 
transformation of ${\mathcal L}^2(\S^1 \times \S^1)$ to the orthocomplement
of the subspace generated by characteristic functions of horizontal sets, the 
spectrum is continuous (see \cite{A}). 
Even though this phenomenon, if occurs at all in our system, can 
only occur for a zero measure set of parameters 
(in the non-Diophantine case), it 
seems to be an interesting point for further investigation.
\end{remark}

\section{Two particles, one rotating scatterer}\label{sec4}
\setcounter{equation}{0}

In this final section we will show that the motion of two particles with one
rotating scatterer, although more complicated, resembles the motion of one
particle with two rotating scatterers.
Let the two point particle velocity vectors 
be $\v$ and $\u$. Note that as we assume 
that only 
the inner scatterer rotates, the outer interactions are billiard collisions
and thus for both 
particles the argument of Lemma~\ref{lemma1}, statement (1) applies.
This means that both of the particles either alternate (all odd bounces are in 
and all even ones are out) or collide always with 
the outer billiard scatterer. We may assume that both particles alternate. 
Otherwise motion is a direct product of two integrable systems: two circle 
billiards if both particles remain outside and a product of a circle billiard 
with the system discussed in Section~\ref{sec2} if one of them alternates.

Recall Remark~\ref{notthesame}: in the present section splitting of $\v$ ($\u$)
into tangential and normal components $(v_t,v_n)$ ($(u_t,u_n)$)
is always understood at the {\it inner} scatterer. 
Again by the geometrical picture
discussed in the proof of Lemma~\ref{lemma1} both $u_n$ and 
$v_n$ are integrals of motion. For brevity we drop subscripts and 
denote $v_t$ ($u_t$) by $v$ ($u$).

We need to introduce one more velocity coordinate $\omega$ for the angular 
velocity of the inner scatterer. Furthermore, as in 
Section~\ref{sec2}, $\eta=\frac{\Theta}{mR^2}$ where 
$\Theta$ is the moment of inertia of the scatterer and $m$ is the mass
of one of the point particles (we assume them to have the same mass).

An interaction between the particle with velocity $\v$ and the inner scatterer
will be referred to as a v-collision, and the definition of u-collision is analogous. 
In course of a v-collision
the velocities transform as
\bea
\label{vcoll}
v' &=& v -\frac{2\eta}{1+\eta}(v-R\omega); \nn \\
R\omega'&=& R\omega + \frac{2}{1+\eta} (v-R\omega).
\eea
We have
\be
\label{vinv}
v+\eta R \omega = v'+ \eta R \omega'; \qquad 
v^2+\eta R^2 \omega^2 = v'^2+  \eta R^2\omega'^2
\ee
while $u$ does not change. For u-collisions  
\bea
\label{ucoll}
u' &=& u -\frac{2\eta}{1+\eta}(u-R\omega); \nn \\
R\omega'&=& R\omega + \frac{2}{1+\eta} (u-R\omega),
\eea
and we have conserved quantities:
\be
\label{uinv}
u+\eta R \omega = u'+ \eta R \omega'; \qquad 
u^2+\eta R^2 \omega^2 = u'^2+  \eta R^2\omega'^2
\ee
in addition to the unchanged $v$. 

Let us introduce $q=\sqrt{\eta}R\omega$. Formulas \eqref{vinv} and 
\eqref{uinv} imply the presence of two integrals:
\be
v^2+u^2+q^2=E \qquad {\rm and} \qquad v+u+\sqrt{\eta}q=N.\nn
\ee 
Thus velocity motion is restricted to the intersection 
of a sphere and a plane in $\R^3$, a circle which will be referred to as the 
velocity circle. The v- and u-collisions, i.e.~the transformations
\eqref{vcoll} and \eqref{ucoll} are reflections across the planes
$$
\sqrt{\eta}v-q=0 \qquad {\rm and} \qquad \sqrt{\eta}u-q=0,
$$
respectively. When restricted to the velocity circle, v-collisions reduce
to reflections across the line
\be
\label{vline}
\big{(} \, N+t\, ,\, -(1+\eta)t\, ,\, \sqrt{\eta}t \, \big{)}
\ee
while u-collisions to reflections across 
\be
\big{(} \, -(1+\eta)t\, ,\, N+t\, ,\, \sqrt{\eta}t \, \big{)}.\nn
\ee
The angle of the two lines, which we denote by $\frac{\gamma}{2}$, is an 
important parameter depending on which we distinguish rational/irrational
cases. We may calculate
\be
\cos\frac{\gamma}{2}=\frac{1}{1+\eta},\nn
\ee
thus, just like in Lemma~\ref{baseangle}, for a full measure set of $\eta$ 
parameters -- actually, apart from countably many points -- the base is 
irrational, and this classification does not depend on the integrals of 
motion.

We  parameterize the velocity circle with the normalized arclength parameter $s\in [0,1)$
as in Convention~\ref{01}.
For definiteness
$s=0$ corresponds to one of the intersection points of the circle with 
the line \eqref{vline}. In this coordinate the v- and u-collisions, 
i.e.~the transformations \eqref{vcoll} and \eqref{ucoll} are
\bea
s &\to& s_1 = -s \  ({\rm mod} \, 1) {\rm \  and}\label{vs}\\ 
s &\to& s_2 = -s-\gamma \ ({\rm mod} \, 1) \label{us},
\eea  
respectively. 

However, at this point the exposition ceases to be parallel to that of
Section~\ref{secvel}. The two transformations above may follow each other
in a more complicated manner and it is not enough to take two copies of the 
circle to describe the velocity evolution.

Consider the Poincar\'e section for the cross section of v-collisions.
When the system
is in state $s$, there is a unique tangential velocity $v(s)$ -- and thus, 
as the 
normal velocity component $v_n$ is an integral of motion, 
a unique velocity vector $\v(s)$ -- for the first particle. Thus the return 
time for the first particle can be calculated (cf. Fig~\ref{fig2}) as
\be
\label{t1}
\tau_1(s)=\tau(v(s)), \quad \hbox{where} \quad \tau(v)=\, 2\, \frac{d(\v)}{|\v|}
\ee
and $d(\v)$ is given in Formula \eqref{Ts} 
(note this time $v_n$ does not depend on $s$ as it is an integral of motion).
The evolution of the system strongly depends on the number of inner bounces the
second particle has between two consecutive bounces of the first one. 
As within this time interval the second particle follows the pattern described 
in Section~\ref{sec2}, we have the following picture. 
Let us assume the system is in state $s\in\S^1$ just after a v-collision and 
thus the angular velocity of the scatterer is $\omega(s)$. 
If the number of u-collisions before the next v-collision is even, 
the scatterer, just before the next v-collision,    
will rotate with the exact same angular velocity $\omega(s)$. 
This implies that 
just after the second v-collision the system arrives at the state
$s_1=-s$ (cf. \eqref{vs}). On the other hand, if the number of 
u-collisions in this time interval is odd, the scatterer rotates with 
angular velocity $\omega(s_2)$ just before the second 
v-collision, and thus the system lands just after the second v-collision
at the state $-s_2=s+\gamma$ (cf. \eqref{us}).

To calculate the number of u-collisions within the time interval $[0,\tau_1]$
we denote by
\be
\label{t2}
\tau_2(s)=\tau(u(s)) \quad \hbox{and} \quad \tau_3(s)=\tau(u(s_2))
\ee
the length of time intervals between consecutive u-collisions in the two 
possible phases of the second particle (cf. also \eqref{t1}). We introduce 
\be
\label{hat}
\hl=\frac{\tau_2}{\tau_2+\tau_3} \ \in [0,1) \quad {\rm and} \quad
\hat{t}=\left\{ \frac{\tau_1}{\tau_2+\tau_3} \right\} \ \in [0,1)
\ee
(where $\{x\}$ is the 
fractional part of $x\in \R$).
If the two particles started from the 
inner scatterer at the same time, the parity of u-collisions until the next 
v-collision would depend only on the relation of these two quantities: 
it would be even/odd for 
$\hat{t}$ less than (greater than) $\hl$. To track 
the velocity evolution in the general case we need to introduce one more 
coordinate. 

\begin{definition}
At the moment of a v-collision (i.e.~on our Poincar\'e section) we define
$t\in[0,1)$ as the following quantity: the last u-collision preceded
the v-collision we are considering at time $t\tau_2(s)$.   
\end{definition}

The coordinate $t$ defined this way describes where the second particle is at 
the moment the first particle collides inside. 

\begin{convention}
\label{deft}
We think of the quantities $t,\hat{t}$ as elements of $\S^1$ (identified with
$[0,1)$) thus we may add them (which is always understood (mod $1$)). 
\end{convention}

This way the velocity evolution can be described as the following two dimensional
system
\bea
T :\S^1\times \S^1 \to \S^1\times \S^1, \quad T(s,t)=(s',t') &\ & 
{\rm where} \nn \\
\label{twobase}
(s',t')=\left( -s, t+\frac{\hat{t}}{\hl} \right) \quad {\rm if} \quad 
t\hl+\t<\hl &\ & {\rm (even \ case)} \\ \nn
(s',t')=\left( s+\gamma \, , \, \frac{\tau_2}{\tau_3}t + 
\frac{\t}{1-\hl} \right) 
\quad {\rm if} \quad 
t\hl+\t>\hl &\ & {\rm (odd \ case)}.
\eea
Note that $T$ is well defined 
since $\t$ and $\hl$ both depend on $s$, see Equation~\eqref{hat}.

\begin{remark}
The reason for the rescaling term $\frac{\tau_2(s)}{\tau_3(s)}$ in 
the evolution of the $t$--coordinate for the odd case is that, as there has
been an odd number of 
u-collisions, the length of the period for the second particle 
is $\tau_3$ in contrast to the previous $\tau_2$.

We have
 $t'=\frac{t\tau_2(s)+\t(\tau_2(s)+\tau_3(s))}{\tau_2(s')}$
for both the odd and the 
even cases. Note $u(s)=u(-s)$ (as the reflection \eqref{vs} corresponds to 
the interaction of the first particle with the scatterer, thus the velocity 
of the second particle does not change) and thus $\tau_2(s)=\tau_2(-s)$ 
(cf. \eqref{t1} and \eqref{t2}). This implies $\tau_2(-s)=\tau_2(s)$ 
(even case) and $\tau_2(s+\gamma)=\tau_3(s)$ (odd case).
\end{remark}

Note that this way we need to study a two dimensional system for 
velocity evolution, which is much more complicated than the one dimensional
base system of Section~\ref{secvel}.
However, as the theorem below shows, for a suitable 
choice of integrals and, correspondingly, initial conditions 
(both of positive measure), there is an arbitrary long time interval
for which the $t$ coordinate can be disregarded and for which the velocity 
motion follows a rotation orbit.    

\begin{theorem}
\label{two}
For any $\e>0$ small enough there exists a positive integer 
$N_{\e}$, an open set $W$ in the space of 
parameters (i.e.~integrals of motion) and a set of initial conditions 
$E$ of Lebesgue measure bigger than 
$\e$ with the following property.
Given a system with parameters from $W$, for any 
$(s_0,t_0)\in E$ the iterate $(s_k,t_k)=T^k(s_0,t_0)$ is 
in the odd part of the phase space, and thus, $s_k=s_0+k\gamma$, as long as
$0\le k<N_{\e}$. Moreover, $N_{\e}\to \infty$ as $\e \to 0$.
\end{theorem}

\begin{proof}
Let us define $E\subset \S^1\times \S^1$ with 
$|t_0-\frac{1}{2}|<\e$ and $s_0\in \S^1$ arbitrary. This means that the second 
particle is close to the outer (billiard) scatterer when the first 
particle collides with the inner (rotating) one.

To define the open set 
$W$ let us fix $N$ and $E$ arbitrary (in such a way that the 
velocity circle is not empty). Note that the normal velocity components are 
integrals of motion as well. We fix $\e'>0$ small 
enough and $K>0$ big enough, their values will be chosen later depending 
on $\e$. For the points of $W$ we assume $|v_n-K'|<\e'$ 
and $|u_n-K'|<\e'$ for some $K'>K$. This ensures that 
the normal velocity components of the particles are almost equal and,
 moreover, much greater than the tangential components. As a consequence,
for all $s$, both $d(v(s))$ and $d(u(s))$ are very close to $1-R$, while 
both $|\v(s)|$ and $|\u(s)|$ are very close to $K'$. In particular,
with $\e'$ sufficiently small and $K$ sufficiently big we may ensure that
\bea
\label{epsilon}
\t(s)\in\left( \frac{1}{2}-\e,\frac{1}{2}+\e \right), \quad &\ & 
\hl(s)\in\left( \frac{1}{2}-\e,\frac{1}{2}+\e \right), \\ \nn
\frac{\tau_2(s)}{\tau_3(s)} \in \left( 1-\e,1+\e \right) , 
\qquad {\rm and} \quad  
&\ & \left| \frac{\t(s)}{1-\hl(s)} \right|< \e
\eea     
for all $s\in\S^1$, cf. \eqref{hat}. Note that in the last 
inequality the absolute value is understood as distance from $0$ on 
$\S^1$.

Finally, we fix 
\be
\label{Neps}
N_{\e}=\left[ \min \left( \, \frac{\log(2-2\e)}{\log(1+3\e)} \, , \, 
\frac{\log(8\e)}{\log(1-3\e)} \, \right) \right] \, - \, 2,
\ee
where $[x]$ denotes the integer part of $x\in \R$.

\begin{lemma}\label{thelemma}
Assume $0\le k <N_{\e}$. We have
$$
(1) \quad \frac{1}{2}(1-3\e)^{k+1}<t_k<\frac{1}{2}(1+3\e)^{k+1}, \qquad
(2) \quad \hl< t_k\hl+\t<1.
$$
Note that (2) means in particular that $(s_k,t_k)$ is in the odd part of
the phase space, cf. \eqref{twobase}.
\end{lemma}
\begin{proof}
We will prove (1) and (2) by induction. As $|t_0-\frac{1}{2}|<\e$, both 
statements hold for $k=0$ by \eqref{epsilon}. 
Assume that we have already established (1) and 
(2) for $k-1$. Note that this, in particular, means that 
$(s_{k-1},t_{k-1})$ is in the odd part and thus by \eqref{twobase}
$$
t_k=\frac{\tau_2}{\tau_3} t_{k-1} + \frac{\t}{1-\hl}\, .
$$
This implies, by the inductive assumption and \eqref{epsilon} 
\bea
t_k <\frac{1}{2}(1+3\e)^k(1+\e)+\e &<& \frac{1}{2} (1+3\e)^{k+1}, \nn \\
t_k >\frac{1}{2}(1-3\e)^k(1-\e)-\e &>& \frac{1}{2} (1-3\e)^{k+1} \nn
\eea
which is (1) for $t_k$. On the other hand, (1) implies (2) as long as 
$k<N_\e$. To see this note that by \eqref{epsilon} it is enough to prove that
\bea
\frac{1}{2}(1-3\e)^{k+1}\frac{1}{2}(1-2\e) +\frac{1}{2}-\e &>& 
\frac{1}{2}+\e \quad {\rm and } \nn \\ \nn
\frac{1}{2}(1+3\e)^{k+1}\frac{1}{2}(1+2\e) +\frac{1}{2}+\e &<& 1.
\eea
These conditions are both satisfied as long as $k<N_{\e}$, cf. 
\eqref{Neps}.
\end{proof}
By Formula~\eqref{twobase}, Theorem~\ref{two} is a straightforward consequence of  
Lemma~\ref{thelemma}.
\end{proof}

\begin{remark}
\label{twoskew}
Considering motion in the fibers, i.e.~in addition to the velocity coordinate 
$s$, the consecutive points of impact of v-collisions (points of impact for
the first particle) parametrized by the arclength parameter $\phi$ along the 
inner scatterer, we can make the following observation. Under the conditions of
Theorem~\ref{two}, the coordinate $t$ can be disregarded, and the two 
dimensional system of the $(s,\phi)$ variables follows (up to time $N_{\e}$)
the Anzai skew product dynamics analyzed in Section~\ref{secskew}. 
\end{remark}

\section*{Acknowledgements}

Useful discussions with B.~Fernandez are thankfully acknowledged. We thank 
S.~Ferenczi for pointing out \cite{A}.   
P.B. is grateful for the hospitality and the inspiring research atmosphere of 
CPT Marseille. This research has been partially supported by the 
Hungarian National Foundation for Scientific Research (OTKA), 
grants T32022 and TS040719.

\end{document}